\documentclass[letter]{amsart}
\usepackage{amsmath,amssymb,amsfonts,amsthm,calc,color}

\usepackage[latin1]{inputenc}

\usepackage[dvipsnames]{xcolor}

\usepackage[german, english]{babel}
\usepackage{hyperref}
\usepackage{enumitem}

\newcommand{\be}{\begin}
\newcommand{\e}{\end}
\newcommand{\beq}{\begin{equation}}
\newcommand{\eeq}{\end{equation}}

\newcommand{\om}{\omega}

\theoremstyle{definition}

\numberwithin{equation}{section}
\newcommand{\comment}[1]{}

\newcommand{\curly}[1]{\mathcal{#1}}
\newcommand{\setof}[2]{\left\{ #1\; : \;#2 \right\}}

\newcommand{\Z}{{\mathbb Z}}
\newcommand{\R}{{\mathbb R}}
\newcommand{\C}{{\mathbb C}}

\newcommand{\N}{{\mathbb N}}

\newcommand{\OO}{{\mathcal{O}}}

\newcommand{\Lp}{{\Delta}}

\newcommand{\al}{{\alpha}}
\newcommand{\de}{{\delta}}
\newcommand{\eps}{{\varepsilon}}

\renewcommand{\l}{\left}
\renewcommand{\r}{\right}

\newcommand{\scp}[2]{\langle#1,#2\rangle}
\newcommand{\scpp}[2]{\l\langle#1,#2\r\rangle}
\newcommand{\qmexp}[1]{\left\langle #1\right\rangle}

\newcommand{\Hm}[1]{\leavevmode{\marginpar{\tiny%
$\hbox to 0mm{\hspace*{-0.5mm}$\leftarrow$\hss}%
\vcenter{\vrule depth 0.1mm height 0.1mm width \the\marginparwidth}%
\hbox to 0mm{\hss$\rightarrow$\hspace*{-0.5mm}}$\\\relax\raggedright
#1}}}

\begin{document}

\title[Asymptotic expansion of annealed Green's function and derivatives]{Asymptotic expansion of the annealed Green's function and its derivatives}

\begin{abstract} We consider random elliptic equations in dimension $d\geq 3$ at small ellipticity contrast. We derive the large-distance asymptotic expansion of the annealed Green's function up to order $4$ in $d=3$ and up to order $d+2$ for $d\geq 4$. We also derive  asymptotic expansions of its derivatives. The obtained precision lies far beyond what is established in 
prior results in stochastic homogenization theory. Our proof builds on a recent breakthrough in perturbative stochastic homogenization by Bourgain in a refined version shown by Kim and the second author, and on Fourier-analytic techniques of Uchiyama.
\end{abstract}

\author[M. Keller]{Matthias Keller}
\address{Matthias Keller, Universit\"at Potsdam, Institut f\"ur Mathematik, 14476  Potsdam, Germany}
\email{matthias.keller@uni-potsdam.de}
\author[M. Lemm]{Marius Lemm}
\address{Marius Lemm, Institute of Mathematics, EPFL, 1015 Lausanne, Switzerland}
\email{marius.lemm@epfl.ch}
\date{July 24, 2021}
\maketitle

\section{Introduction}
The regularity of Green's functions and their derivatives forms the backbone of classical elliptic regularity theory for divergence-form operators \cite{deGiorgi,LSW63,Moser,Nash}. Here we consider divergence-form discrete elliptic operators of the form
\beq\label{eq:ellipticop}
\nabla^* \mathbf A_\omega\nabla, \qquad \textnormal{on } \ell^2(\Z^d), \, d\geq 3,
\eeq
with a random, elliptic coefficient matrix $\mathbf A_\omega(x)\in \R^{d\times d}$. The central goal of the very lively field of stochastic homogenization theory is to understand the \textit{large-distance behavior} of solutions to $\nabla^* \mathbf A_\omega(x) \nabla u_\omega(x) =f(x)$. Naturally, Green's functions play a central role in this endeavour; see \cite{AKMpaper,AKM,BGO1,GNO} and references therein.

In 2018, J.~Bourgain \cite{B} introduced a completely novel approach to studying such equations in the regime of small ellipticity contrast. He takes the coefficients to be
\beq\label{eq:coefficients}
 \mathbf A_\omega(x)=(1+\de \sigma_\om(x))\mathbf I_d
\eeq
with $\{\sigma_\om(x)\}_{x\in\Z^d}$ a family of independent and identically distributed bounded random variables, $\mathbf I_d$ the $d\times d$ identity matrix, and $\de>0$ is a small parameter. A key point is that inspired by an earlier unpublished note of I.M.~Sigal, from the outset Bourgain's focus lies not with deriving an effective large-distance description of (random) solutions $u_\om(x)$, but only of their \textit{average} $\qmexp{u_\om(x)}$. This is equivalent to studying the annealed (i.e., averaged) Green's function
$$
\mathcal G(x)=\qmexp{G_\omega(x,0)}=\qmexp{\frac{1}{\nabla^* \mathbf A_\omega \nabla}\de_0(x,0)},
$$
where $ 1/ \nabla^* \mathbf A_\omega \nabla$ is the operator inverse of $ \nabla^* \mathbf A_\omega \nabla $, see  \cite{KL}.
The main result of \cite{B}, which was subsequently refined by Kim and the second author \cite{KL}, establishes that $\qmexp{G_\omega(x,0)}$ can be represented as a convergent perturbation series in $\de>0$ with explicit large-distance decay bounds. (See Theorem \ref{thm:bkl} below for the precise statement.) This fact has several non-trivial consequences. For instance, \cite{DGL} showed that it allows to define higher-order correctors beyond what was previously believed possible. A related, intriguing possibility is to extend the main result of \cite{B,KL} (see Theorem~\ref{thm:bkl} below) to arbitrary ellipticity contrast. This is known as the Bourgain-Spencer conjecture which remains open. See \cite{D,DLP} for recent results in this direction.

For our purposes here, we stay within the small-ellipticity contrast regime and instead focus on a different consequence of the main results in \cite{B,KL}. Namely, \cite[Corollary 3.1]{KL} proves that for all multi-indices $\al$ with $|\al|\leq d+1$,
\beq\label{eq:Galphaestimate}
|\nabla^\al \mathcal G(x)|\leq C_\alpha (1+|x|)^{2-d-|\al|},\qquad x\in \Z^d.
\eeq
That is, the first $d+1$ derivatives match the expected power-law scaling that is familiar from the free Laplacian. (Before these works, \eqref{eq:Galphaestimate} was only known for $|\al|\leq 2$, in any dimension and for any ellipticity contrast \cite{AKM,DD,MO1,MO2}.)\medskip

\textit{In the present work, we build further on the results of \cite{B,KL} to prove a previously unforeseen strong refinement of \eqref{eq:Galphaestimate}, a precise asymptotic expansion of $\nabla^\al \mathcal G(x)$ as $|x|\to\infty$.}
\medskip
 
Our main results can be summarized as follows.

\be{itemize}
\item Theorem \ref{thm:GFasymptotic} provides the asymptotic expansion of $\mathcal G(x)$ as $|x|\to\infty$ up to order $4$ in $d=3$ and up to order $d+2$ for $d\geq 4$. 
\item Corollary~\ref{cor:nablaGFasymptotic} contains analogous asymptotic expansions for the derivatives $\nabla^\alpha \mathcal G$ with $|\al|\leq 3$ in $d=3$ and $|\al|\leq d+1$ in $d\geq 4$. For every derivative taken, one loses one order in the asymptotic expansion of $\mathcal G(x)$, so when one reaches the last derivative, we only identify the leading-order asymptotics.
\e{itemize}

In general, these main results go far beyond what can be achieved by the powerful methods of homogenization theory in the regime of small ellipticity contrast. The reason is partly that those methods first describe the random Green's function $G_\om(x,0)$ which is harder to understand due to probabilistic fluctuations, see e.g.\ \cite[Theorem 8.20 \& Section 9.2]{AKM}, \cite[Corollary 3]{BGO1}, \cite[Proposition 4.2]{GM} and \cite[Theorem 5.1]{MO3}. These results can be averaged post-hoc to obtain information on $\mathcal G(x)=\qmexp{G(x,0)}$ and its derivatives. For instance, it is well-known that the leading term in the expansion is the homogenized Green's function; see \cite{AKM,MO1} and the other references above. However, the resulting bounds will be much less precise than the expansion we show here in Theorem \ref{thm:GFasymptotic} or Corollary~\ref{cor:nablaGFasymptotic}. In particular, there does not appear to be any result on derivatives of order $\geq 3$ in the literature.

The paper is organized as follows.

\be{itemize}
\item In \textbf{Section \ref{sect:main}}, we define the setting and state the main results.
\item In \textbf{Section \ref{sect:mainthmpf}}, we prove Theorem \ref{thm:GFasymptotic}.
\item In \textbf{Section \ref{sect:maincorpf}}, we prove Corollary~\ref{cor:nablaGFasymptotic}.
\item  In the \textbf{Appendix}, we include a self-contained derivation of the leading term in Theorem \ref{thm:GFasymptotic} based on \cite[Theorem 1.1]{KL}.
\e{itemize}

An open question related to diffusion processes is discussed in Subsection \ref{ssect:T}.

We mention that a version of the main results of this paper originally appeared in the preprint \cite[Version 1]{KLhardy} about optimal Hardy weights on $\Z^d$. The present paper has been split off from that work.

\section{Setup and main results}\label{sect:main}
We begin by reviewing the setup and the main results of \cite{B,KL} which form the backbone of our asymptotic expansion. Afterwards, we state our main results, Theorem \ref{thm:GFasymptotic} and Corollary~\ref{cor:nablaGFasymptotic}.

\subsection{Basic setting}
 Recall the definition of the discrete derivative denoted by $\nabla=(\nabla_1, \nabla_2,\ldots,\nabla_d)^T$. For a function $u:\Z^d\to \R$ or $\C$,
$$
\nabla_j u (x) = u(x+e_j)-u(x)
$$
with $e_j$ the $j$th canonical basis vector. Its $\ell^2(\Z^d)$-adjoint is denoted by $\nabla^*=(\nabla_1^*, \ldots,\nabla_d^*)$ and acts as
$$
\nabla_j^* u (x) := u(x-e_j)-u(x).
$$
Then $\nabla^*\nabla=-\Delta$ is the usual discrete Laplacian, a positive operator.

\be{assumption}\label{ass:iid}
Let $\{\sigma_\om(x)\}_{x\in\Z^d}$ be a family of independent and identically random variables bounded by $1$.
\e{assumption}

We shall consider the random divergence-form operator 
\beq\label{eq:operator}
L_\om=-\Delta + \delta \nabla^* \sigma_\om \mathbf I_d \nabla, \qquad \textnormal{on } \ell^2(\Z^d).
\eeq
We will suppress the identity matrix $\mathbf I_d$ from the notation. Note that $L_\om$ is of the form \eqref{eq:ellipticop} with the coefficients chosen by \eqref{eq:coefficients}. We assume that $\de\in (0,1)$ so that $L_\om$ is uniformly elliptic. We denote the Green's function of $L_\om$ by $G_\om(x,0)$, the unique solution to $L_\om G_\om (x,0)=\de_0(x)$ which is well-defined for $d\geq 3$. 
 
 Our main object of interest is the \textit{annealed Green's function}
$$
\mathcal G(x)=\qmexp{G(x,0)}.
$$
The relevance of the annealed Green's function is that it governs the behavior of averaged solutions. For instance, take \cite[Corollary 1.6]{KL}. It says that for any  $f\in\ell^{p_d}(\Z^d)$ with $p_d^{-1}=\tfrac{1}{2}+\tfrac{1}{d}$ (the critical Sobolev index), there exists a unique random solution $u_\om \in \ell^{q_d}(\Z^d)$, with  the H\"older dual $q_d$ of $p_d$ so that $u_\om$ solves the equation $L_\om u_\om=f$ and the averaged solution is given by the formula
$$
\qmexp{u_\om}=\mathcal G*f.
$$
Thus we see that \textit{the decay properties of $\mathcal G$ determine the decay properties of the averaged solution $\qmexp{u_\om}$}.  The same is true for derivatives of all orders. 

The correspondence is cleanest when $f$ is compactly supported. Note that decay rates of a function and its derivatives are the most natural way to measure regularity on $\Z^d$. 

\be{remark}
The coefficients $\sigma_\om$ are taken to be a multiple of the identity matrix $\mathbf I_d$ only for simplicity. The same techniques apply if the  i.i.d.\ perturbation is any symmetric matrix \cite[Remark~1.4]{KL}. 
\e{remark}

\subsection{Background on the annealed Green's function}
In \cite{B}, Bourgain shows that the annealed Green's function arises itself as a Green's function of a matrix-valued convolution operator, called $\curly{L}$ below, which arises as the harmonic mean of the original random operator. This ``parent operator'' for $\mathcal G$ can be realized as a bounded operator 
$$
\mathcal L:H^1(\Z^d)\to H^{-1}(\Z^d)
$$
defined via the discrete Sobolev spaces 
$$
H^1(\Z^d)=\Lambda^{-1}(\ell^2(\Z^d)),\qquad H^{-1}(\Z^d)=\Lambda(\ell^2(\Z^d)),\qquad \Lambda =(-\Delta)^{1/2}.
$$
We refer to \cite[Section 2.1]{KL} for the details and to \cite[Lemma 1.1]{DGL} for an alternative definition of $\mathcal L$ via the Lax-Milgram theorem.

The breakthrough result of \cite{B} gives a precise description of $\mathcal L$ of the following form
\beq\label{eq:harmonicmeanoperator}
\curly{L}=\Lp +\nabla^* \mathbf{K}^\de \nabla,
\eeq
where $ \Delta $ is the free Laplacian and  $\mathbf{K}^\de$ is a  $d\times d$ matrix-valued convolution operator whose components satisfy a decay estimate. This decay estimate was subsequently improved to the (conjecturally nearly optimal) rate $-3d+\eps$ in \cite{KL} which we use here.

We now summarize these results. We notationally identify the convolution operator $\mathbf{K}^\de$ with its matrix-valued kernel $\mathbf{K}^\de(x-y)\in \R^{d\times d}$.

\be{thm}[\cite{B,KL}]\label{thm:bkl}
Let $d\geq 3$ and $\eps\in (0,1)$. There exists $c_d>0$ so that for all $\de\in (0,c_d \eps)$, the representation \eqref{eq:harmonicmeanoperator} holds with the following decay estimate on the convolution kernel
$$
|K^\de_{j,k}(x-y)|\leq C_{d} \de^2 (1+|x-y|)^{-3d+\eps},\qquad  j,k\in \{1,\ldots,d\}.
$$
\e{thm}

For the purposes of this paper, we can choose $\eps=\tfrac{1}{2}$ for definiteness so that, for $\de\in (0,c_d)$,
\beq\label{eq:Kdeltadecay}
|K^\de_{j,k}(x-y)|\leq C_{d} \de^2 (1+|x-y|)^{-3d+\eps},\qquad  j,k\in \{1,\ldots,d\}.
\eeq

%
%For $d\geq 3$, there exists a unique self-adjoint convolution operator $B:\ell^2(\Z^d)\to\ell^2(\Z^d)$ with $1-\de\leq B\leq 1+\de$ that satisfies the following. For all vector-valued functions $f\in \ell^2(\Z^d)^d$, if
%$$
%(-\Delta + \delta \nabla^* \sigma_\om \nabla) u_\om(x)=\nabla\cdot  f(x),
%$$
%then the average $\qmexp{u_\om(x)}\in \dot{H}^1(\Z^d)$ is the unique Lax-Milgram solution of
%$$
%B\qmexp{u_\om(x)} =f(x),
%$$
%see \cite[Lemma 1.1]{DGL} and also \cite[Section 2.1]{KL}.

The usefulness of \eqref{eq:Kdeltadecay} lies in the fact that it guarantees the existence of moments of $\mathbf{K}^\delta$ up to order $2d-1$.
The existence of higher moments has meaning in homogenization theory, where it can be shown to be equivalent to the existence of a previously unforeseen higher-order corrector theory up to order $2d$, \cite{DGL}.

Of particular importance for the leading-order behavior is the $d\times d$ matrix
\beq\label{eq:Qdefn}
\mathbf Q=\mathbf{I}_d+\sum_{x\in\Z^d}\mathbf{K}^\delta(x).
\eeq
In the language of homogenization theory, 
\beq\label{eq:Qhom}
\mathbf Q=\frac{\overline{\mathbf{a}}+(\overline{\mathbf{a}})^T}{2}
\eeq
corresponds to the symmetrized \textit{lowest-order homogenized coefficients} \cite[Eq.\ (2.5)]{DGL}. %This also shows that $\mathbf{Q}$ symmetric.

\be{proposition}\label{prop:Qsymm}
The $d\times d$ matrix $\mathbf{Q}$ is symmetric
 and there exist constants $c_d,C_d>0$ so that for all $\de\in (0,c_d)$,
\beq\label{eq:Qbounds}
1-C_d \de^2\leq \mathbf{Q}\leq 1 +C_d\de^2.
\eeq
\e{proposition}

\be{proof}
For the symmetry of  $\mathbf{Q}$, we use the power series representation of $\mathbf{K}^\delta$, cf.\ \cite[Eq.~(2.5)]{B} and \cite[Eq.~(1.14)]{KL}. To write this down, we require some basic objects and notation from \cite{B,KL}. Let $\Omega$ denote the underlying probability space to the $\{\sigma_\om(x)\}_{x\in Z^d}$. 

First, we express the expectation as a projection operator on the extended space
$$
\begin{aligned}
P:\;L^2(\Z^d\times \Omega)\to& \ell^2(\Z^d)\subset L^2(\Z^d\times \Omega)\\
u(x,\omega)\mapsto& \qmexp{u}(x)
\end{aligned}
$$
Here $L^2(\Z^d\times \Omega)$ is defined with respect to the counting measure on $\Z^d$ and the probability measure on $\Omega$. We write $P^\perp=\mathbf{I}_{L^2(\Z^d\times \Omega)}-P$ for the projection onto the orthogonal complement.

Second, we write $\sigma$ for the multiplication operator 
$$\sigma :L^2(\Z^d\times \Omega)\to L^2(\Z^d\times \Omega) ,\quad (\sigma u)(x,\om)= \sigma_\om(x) u(x,\om).$$

Third, we introduce the operator-valued $d\times d$ matrix $\mathbf{K}$ whose components are the operators $\mathbf{K}_{j,k}:\ell^2(\Z^d)\to\ell^2(\Z^d)$ which are defined as Fourier multiplication by the functions 
$$
F(\theta)=\frac{(e^{i\theta_j}-1)(e^{-i\theta_k}-1)}{2\sum_{j=1}^d(\cos\theta_j-1)}
$$
using the following convention for the Fourier transform
$$
\hat f(\theta)= \sum_{x\in\Z^d} e^{-ix\cdot\theta}f(x),\qquad \theta\in [-\pi,\pi]^d.
$$
(The operator $\mathbf{K}$ can be formally written as $\frac{\nabla\nabla^*}{\Delta}$ and is also known as the discrete Helmholtz projection.)  Equivalently, the operator $\mathbf{K}_{j,k}$ is a convolution operator with the convolution kernel 
\beq\label{eq:Kkernel}
\mathbf{K}_{j,k}(x-y)=\int_{[-\pi,\pi]} e^{i(x-y)\cdot\theta}F(\theta)\frac{\mathrm{d}^{d}\theta}{(2\pi)^d}.
\eeq
Then we lift $\mathbf{K}$ to the extended space $L^2(\Z^d\times \Omega)$ by acting trivially on the random component and we abuse notation by calling the resulting operator $\mathbf K$ as well.

With these preparations complete, we can write down the following power series representation of $\mathbf{K}^\delta$, cf.\ \cite[Eq.~(1.14)]{KL}
\beq\label{eq:pertseries}
\mathbf{K}^\delta=\de \sum_{n=1}^\infty (-\de)^n P\sigma (\mathbf{K} P^\perp \sigma)^n.
 \eeq
 We remark that the convergence of this series is proved in \cite{B,KL}. 
 
 We claim that \eqref{eq:pertseries} implies that the concolution kernel satisfies
 \beq\label{eq:KdeltaT}
( \mathbf{K}^\delta)^T(x)=\mathbf{K}^\delta(-x).
 \eeq
Let us prove this. By \eqref{eq:pertseries}, we have
$$
\begin{aligned}
&( \mathbf{K}^\delta)^T(x)\\
=&\de \sum_{n=1}^\infty (-\de)^n P\sigma (\mathbf{K}^T P^\perp \sigma)^n(x,0)\\
=&\de \sum_{n=1}^\infty (-\de)^n P\sigma(x)
\sum_{x_1,\ldots,x_{n-1}\in\Z^d}
\mathbf{K}^T(x-x_1) P^\perp \sigma(x_1)
\ldots 
\mathbf{K}^T(x_{n-1}) P^\perp \sigma(0).
\end{aligned}
$$
By \eqref{eq:Kkernel}, we have $\mathbf{K}^T(x)=\mathbf{K}(-x)$.  Introducing the reflected configuration $\tilde\sigma_\om(x)=\sigma_\om(-x)$ and reflecting the summation variables $x_\al$ to   $ -x_\al$, we obtain
\beq\label{eq:tildesigappears}
( \mathbf{K}^\delta)^T(x)=\de \sum_{n=1}^\infty (-\de)^n P\sigma (\mathbf{K} P^\perp \tilde\sigma)^n(-x,0).
\eeq
Thanks to Assumption \ref{ass:iid}, $\tilde\sigma$ has the same distribution as $\sigma$. Due to the presence of the first $P$ projection in \eqref{eq:tildesigappears}, $\tilde\sigma$ only appears in an averaged sense in that equation and so \eqref{eq:KdeltaT} is proved.  

Combining \eqref{eq:Qdefn} and \eqref{eq:KdeltaT} with the change of variables $x\to -x$, we conclude
$$
\mathbf{Q}^T=\mathbf{I}_d+\sum_{x\in \Z^d}\left(\mathbf K^\delta\right)^T(x)=\mathbf{I}_d+\sum_{x\in \Z^d} \mathbf K^\delta(x)=\mathbf{Q}
$$
as desired.

Finally, the decay estimate \eqref{eq:Kdeltadecay} implies $ \|\widehat{\mathbf{K}^\delta}(0)\|\leq C_d\de^2$ and so \eqref{eq:Qbounds} follows from the spectral theorem.
\e{proof}

From now on we assume that $ \delta $ is sufficiently small such that \eqref{eq:Kdeltadecay} holds and $ \mathbf{Q} $ is positive definite. 

\be{remark}
A curious but rather unexplored property of $\mathbf {K}^\delta (x)$ is that, despite being an averaged object, it holds enough information to fully characterize the law of the probablity measure of the random coefficients $\{\om_x\}_{x\in\Z^d}$ \cite[Proposition 1.8]{KL}.
\e{remark}

\subsection{Main result 1: Green's function asymptotics}
The asymptotic expansion involves the modified spatial variable
\beq\label{eq:tildexdefn}
\tilde x = \sigma \mathbf{Q}^{-1/2} x,\qquad \textnormal{with } \sigma=(\det\mathbf{Q})^{1/(2d)},
\eeq
and the universal constant
\beq\label{eq:kappadefn}
\kappa_d=\frac{1}{2}\pi^{-d/2}\Gamma(d/2-1).
\eeq

For $d\geq 3$, we denote
\beq\label{eq:mddefn}
m_d=
\be{cases}
3,\qquad &\textnormal{if } d=3,\\
d+1,\qquad &\textnormal{if } d\geq 4.
\e{cases}
\eeq
We now state our first main result, a large-distance asymptotic expansion of $\mathcal G(x)=\qmexp{G(x,0)}$ of order $m_d+1$.

\be{thm}[Asymptotic expansion of the annealed Green's function]\label{thm:GFasymptotic}
Let $d\geq 3$. There exists $c_d>0$ so that for all $\de\in (0,c_d)$ the following holds. There are polynomials $U_1,\ldots,U_{m_d}$ with $U_k$ having degree at most $3k$ so that
\beq\label{eq:GFasymptotic}
\mathcal G(x)=\frac{\kappa_d}{\sigma^2} |\tilde x|^{2-d}
+\sum_{k=1}^{m_d} U_k\l(\frac{\tilde x}{|\tilde x|}\r) |\tilde x|^{2-d-k} +o(|\tilde x|^{2-d-m_d}), \qquad \textnormal{as } |x|\to\infty.
\eeq
\e{thm}

The proof of Theorem \ref{thm:GFasymptotic} is given in Section \ref{sect:mainthmpf}. 

The polynomials $U_k$ are defined in an explicit manner following \cite{Uch98}. They are given as Fourier transforms of fractions of the form $\tfrac{P_{2d-2+k}(\xi)}{|\xi|^{2d}}$ where $P_{2d-2+k}$ is a homogeneous polynomial of degree $2d-2+k$. These polynomials can be explicitly computed as moments of the function
$T:\Z^d\to\R$ defined by
\beq\label{eq:Tdefn}
\begin{aligned}
T(x)=\frac{1}{2}\de_{x=0}+\frac{1}{4d}\de_{|x|=1}+\frac{1}{4d}\sum_{j,k=1}^d \big(&-K^\de_{j,k}(x)+K^\de_{j,k}(x-e_j)\\
&+K^\de_{j,k}(x-e_k)-K^\de_{j,k}(x-e_j-e_k)\big).
\end{aligned}
\eeq
 For instance, we have
$$
\begin{aligned}
U_1(\omega)=& \int_{\mathbb R^d} \frac{P_{2d-1}(\xi)}{|\xi|^{2d}} e^{-i\omega \cdot \xi}\mathrm{d} \xi,\\
P_{2d-1}(\xi)
=&-\frac{2i}{3\sigma^4 (2\pi)^d} |\xi|^{2d-4} \sum_{x\in\Z^d} T(x) (\xi\cdot x)^3. 
\end{aligned}
$$
Here the Fourier transform giving $U_1(\omega)$ is defined in the sense of tempered distributions on $\mathbb R^d\setminus\{0\}$ and the integral can be computed via \cite[Lemma 2.1]{Uch98}. 

\be{remark} \label{rmk:main1}
\be{enumerate}[label=(\roman*)]
\item In the context of homogenization theory, the leading term in the expansion \eqref{eq:GFasymptotic},
$$
\frac{\kappa_d}{\sigma^2} |\tilde x|^{2-d}=G_{\mathrm{hom}}(x),
$$
is well-known as the \textit{homogenized Green's function} $G_{\mathrm{hom}}(x)$. 
\item By \eqref{eq:tildexdefn} and Proposition \ref{prop:Qsymm}, we have the comparability
\beq\label{eqn:xtilde}
C_{d,\de}^{-1}|x|\leq |\tilde x|\leq C_{d,\de}|x|
\eeq
 for an appropriate constant $C_{d,\de}>1$. In particular, $|x|\to\infty$ and $|\tilde x|\to\infty$ are equivalent and $o(|x|^{-k})=o(|\tilde x|^{-k})$.

\item As described above, the polynomials $U_k$ are computable from moments of the operator $\mathbf{K}^\delta$ alone. In view of the results in \cite{DGL}, it is possible to rephrase the asymptotic expansion in terms of higher-order correctors. This reformulation could pave the way for extending Theorem \ref{thm:GFasymptotic} beyond the small-ellipticity regime, i.e., to all $\de\in (0,1)$, as further progress is made on the Bourgain-Spencer conjecture, cf.\ \cite{D}.

\item In Appendix \ref{sect:direct}, we give a a short self-contained argument for readers interested in seeing how the matrix $\mathbf{Q}$ from Theorem \ref{thm:bkl} arises in $\mathcal G(x)$ by Taylor expansion around the origin in Fourier space. This yields the leading asymptotic order in Theorem \ref{thm:GFasymptotic}. The idea for this is to reduce $\mathcal G(x)$ to the Green's function of the free Laplacian, using dyadic pigeonholing to control the error terms. 
\e{enumerate}
\end{remark}

\subsection{Main result 2: Asymptotics of Green's function derivatives}
%First, we remark that the polynomials $U_1,\ldots,U_{m_d}$ are in principle computable from knowledge of moments of $\mathbf{K}^\delta$ alone, cf.\ \cite{Uch98}. 

In the discrete setting, pointwise asymptotics up to order $N$ of a function yield pointwise asymptotics of its first derivatives up to order $N-1$. This procedure can be iterated for higher derivatives, with a loss of one asymptotic order per derivative. Since our expansion \eqref{eq:GFasymptotic} has $m_d+1$ terms, we can describe asymptotics of the derivatives $\qmexp{\nabla^\alpha G}$ with $|\alpha|\leq m_d$ up to order $m_d-|\alpha|$ with $m_d$ defined in \eqref{eq:mddefn}. To this end, we recall the notation of the discrete derivative $\nabla_j$ 
$$
\nabla_j u(x)= u(x+e_j)-u(x).
$$
For a given multi-index $\alpha=(\alpha_1,\cdots,\alpha_d)\in \Z^d$, $\alpha_j\geq 0$, we write $|\alpha| = \sum_{j=1}^d \alpha_j$ and $$ \nabla^\alpha=\nabla_1^{\alpha_1} \cdots \nabla_d^{\alpha_d}. $$
By linearity, $\nabla^\al \mathcal G(x)=\qmexp{\nabla^\al G(x,0)}$.
%so when we write $\nabla G(b)$ with $b\in \mathbb E^d$ it is implicit that the derivative is taken in the first argument. %For an edge $b\in\mathbb E^d$, we let $|b|$ denote the minimal Euclidean distance between its endpoints and the origin.

\be{corollary}[Asymptotic expansion of derivatives of $\mathcal G$]\label{cor:nablaGFasymptotic}
Under the same assumptions as in Theorem \ref{thm:GFasymptotic}, let $\alpha\in \N_0^d$ be a multi-index with $|\al|\leq m_d$. %\beq\label{eq:nablaGFasymptotic}
%\qmexp{\nabla G([x,x+s e_j])}=s \frac{2-d}{2}\frac{\kappa_d}{\sigma^2} |\tilde x|^{1-d} \frac{\tilde x_j}{|\tilde x|}
%+\OO(|\tilde x|^{-d}), \qquad \textnormal{as } |x|\to\infty.
%\eeq 

Then, as $ |x|\to\infty$,
\beq\label{eq:nablaalphaGFasymptotic}
\begin{aligned}
\nabla^\alpha \mathcal G(x)=&\nabla^\alpha\l(\frac{\kappa_d}{\sigma^2} |\tilde x|^{2-d}
+\sum_{k=1}^{m_d-|\alpha|} U_k\l(\frac{\tilde x}{|\tilde x|}\r) |\tilde x|^{2-d-k}\r)+o(|\tilde x|^{2-d-m_d}).
\end{aligned}
\eeq
\e{corollary}

\be{remark}\
 \label{rmk:main2}
\be{enumerate}[label=(\roman*)]
\item Since each $U_k$ is a polynomial and thus smooth, we can apply the mean value theorem to bound the discrete derivative by the corresponding continuum derivative, cf.\ \cite[Lemma, p.~6]{Don}, taking into account the linear change of variables $\tilde x = \sigma \mathbf{Q}^{-1/2} x$. Together with Proposition~\ref{prop:Qsymm}, this readily implies that the decay rates of each term is controlled by
\beq\label{eq:Uorder}
\nabla^\alpha \l(U_k\l(\frac{\tilde x}{|\tilde x|}\r) |\tilde x|^{2-d-k}\r)=O(|\tilde x|^{2-d-k-|\alpha|}).
\eeq
In summary, \eqref{eq:nablaalphaGFasymptotic} indeed gives an asymptotic expansion comprising $m_d-|\alpha|$ orders. 
\item The leading term does not involve $U_k$ and is therefore particularly easy to compute. For example, we have the following leading-order gradient asymptotic 
\beq\label{eq:nablaGFasymptotic}
\qmexp{\nabla G([x,x+s e_j])}=s \frac{2-d}{2}\frac{\kappa_d}{\sigma^2} |\tilde x|^{1-d} \frac{\scp{\tilde x_j}{\tilde e_j}}{|\tilde x|}
+\OO(|\tilde x|^{-d}), \qquad \textnormal{as } |x|\to\infty.
\eeq 
We give a proof of this fact in Section~\ref{sect:maincorpf}.
\item For the maximal value $|\al|=m_d$, \eqref{eq:nablaalphaGFasymptotic} reduces to
\beq
\qmexp{\nabla^\alpha G([x,x+s e_j])}=\frac{\kappa_d}{\sigma^2} \nabla^\alpha  |\tilde x|^{2-d}+o(|\tilde x|^{2-d-m_d}), \qquad \textnormal{as } |x|\to\infty,
\eeq
so it just manages to capture the leading order asymptotic of the $m_d$-th derivative. 
\e{enumerate}
\e{remark}

\subsection{An open question: Does $\curly{L}$ generate a random walk?}
\label{ssect:T}
We close the presentation of the main results by describing an interesting open problem.

In the setting of \cite{B,KL} described above, one may ask whether the operator $\curly{L}$ from \eqref{eq:harmonicmeanoperator} is again the generator of a random walk on $\Z^d$. More precisely, one can write $m(\theta)=4d(1-\hat T(\theta))$ with the function $T$ defined in \eqref{eq:Tdefn}

An interesting simple-to-state question is then the following: Is $T(x)\geq 0$ for all $x\in\Z^d$? If so, $T(x)$ can be interpreted as the transition function of a random walk with generator $\mathcal L$, at least up to a multiplicative factor $4d$. This would mean that probabilistic averages of solutions are themselves governed by \textit{bona fide} diffusion process, whose dynamics may in turn hold non-trivial information about the non-averaged processes.

Such a direct dynamical meaning of the operator $\mathcal L$ is not at all obvious and would be remarkable. We encountered this question when noting that Uchiyama's analysis  \cite{Uch98} would apply more directly if $T(x)\geq 0$. Our initial investigations indicate that identifying the conditions under which $T(x)\geq 0$ is true is connected to subtle questions concerning componentwise positivity of matrix inverses \cite{JS}.

\section{Proof of Theorem \ref{thm:GFasymptotic}}\label{sect:mainthmpf}
The proof relies on Theorem \ref{thm:bkl}, specifically the decay estimate \eqref{eq:Kdeltadecay} which is the main result of \cite{KL}, and the generalization of a delicate Fourier analysis developed by Uchiyama \cite{Uch98} in the probabilistic setting of random walks. 

\subsection{Fourier-space representation}

Recall that the averaged Green's function $\qmexp{G}$ is the Green's function of the operator $ \mathcal{L} $
which can be described via Theorem~\ref{thm:bkl} as \eqref{eq:harmonicmeanoperator} and the decay estimate \eqref{eq:Kdeltadecay}.

Equivalently, the operator $\curly{L}$ is a Fourier multiplier with the symbol $m: \mathbb T^d\to \mathbb C$, on the torus $\mathbb T^d=(\R/2\pi \Z)^d$, given by
\beq\label{eq:mthetadefn}
m(\theta)=2\sum_{j=1}^d (1-\cos\theta_j)+\sum_{1\leq j,k\leq d}
 (e^{-i\theta_j}-1)\widehat{K^\de_{j,k}}(\theta) (e^{i\theta_k}-1).
\eeq
By integration by parts, the decay bound \eqref{eq:Kdeltadecay} then implies the regularity
 $$ 
 \widehat{K^\de_{j,k}}\in C^{2d-1}(\mathbb T^d),\qquad j,k\in \{1,\ldots,d\}
 $$
which will be used many times in the following argument.

By Taylor expansion of the Fourier multiplier \eqref{eq:mthetadefn} at the origin, we find that the lowest order is quadratic and given by \eqref{eq:Qdefn}, i.e., 
$$
\mathbf Q=\mathrm{Hess}(m)(0)=\mathbf{I}_d+\widehat{\mathbf{K}^\delta}(0).
$$

This also means we can express the averaged Green's function as a Fourier multiplier.
\beq\label{eq:GavgFourierrep}
\mathcal G(x)=\int_{\mathbb T^d} e^{ix\cdot \theta} \frac{1}{m(\theta)} \frac{\mathrm{d}^d \theta}{(2\pi)^d}
\eeq
To see that the integral is well-defined for $d\geq 3$, observe that for $ \delta $ small enough $ m $ vanishes only at the origin $ \theta=0 $ by \eqref{eq:Kdeltadecay}. Thus, by Taylor expansion, $\frac{1}{m(\theta)}$ only has a quadratic singularity at the origin.)

The goal is to perform asymptotic analysis of \eqref{eq:GavgFourierrep} as $|x|\to\infty$. This is a delicate stationary phase argument which has to take special care of the singularity at the origin in Fourier space. A hands-on approach to obtain the leading term which is based on the harmonic analysis ideas in \cite[Appendix~A]{KL} is explored in the Appendix. To derive the full asymptotic expansion, we draw on the techniques of Uchiyama \cite{Uch98} who elegantly accounts for cancelations of naively non-integrable terms. While Uchiyama assumes he is in a probabilistic setting which may not pertain to the averaged Green's function, cf.\ Section \ref{ssect:T}, we show now that his argument extends to our case. %(This is perhaps not so surprising considering that there are only few Fourier-theoretic arguments that leverage positivity conditions.)

 To make contact with the probabilistic perspective, we denote
\beq\label{eq:mrewrite}
m(\theta)=4d\l(1-\hat T(\theta)\r)
\eeq
with $T:\Z^d\to\R$ given as in Section \ref{ssect:T}, i.e.,
\beq\label{eq:Tdefn'} 
\begin{aligned}
T(x)=\frac{1}{2}\de_{x=0}+\frac{1}{4d}\de_{|x|=1}+\frac{1}{4d}\sum_{j,k=1}^d \big(&-K^\de_{j,k}(x)+K^\de_{j,k}(x-e_j)\\
&+K^\de_{j,k}(x-e_k)-K^\de_{j,k}(x-e_j-e_k)\big).
\end{aligned}
\eeq
Note that we produced the term $\frac{1}{2}\de_{x=0}$ by adding and subtracting a constant in \eqref{eq:mrewrite}. This is a common technical trick in the context of discrete random walks to remove periodicity, cf.\ \eqref{eq:ap} below.\medskip

\subsection{Properties of $T$} As mentioned above, our goal is to extend  \cite[Theorem~2]{Uch98} to our situation. In a first step, we verify the assumptions of that theorem with the exception of $T(x)\geq 0$. The function $T$ satisfies the following properties assumed in \cite{Uch98} for small $ \delta $. For all these properties the decay bound  \eqref{eq:Kdeltadecay}, which is $ |K^{\delta}_{j,k}(x)|\leq C_{d}\delta^{2} (1+|x|^{-3d+1/2}) $,  from \cite[Theorem 1.1]{KL} is of the essence. 

\begin{enumerate}[label=(\roman*)]
\item  $T$ has zero mean. Indeed, by \eqref{eq:Kdeltadecay}, we can use Fubini and a change of variables to see
\beq
\label{eq:Uass1}
\begin{aligned}
\sum_{x\in\Z^d} xT(x)
=&\frac{1}{4d}\sum_{x\in\Z^d} x \de_{|x|=1}\\
&+\frac{1}{4d}\sum_{j,k=1}^d \sum_{x\in\Z^d} K^\de_{j,k}(x) \l(-x+(x+e_j)+(x+e_k)-(x+e_j+e_k)\r)\\
=&0.
\end{aligned}
\eeq

\item The smallest subgroup of $\Z^d$ generated by 
\beq\label{eq:ap}
\setof{x\in\Z^d}{T(x)>0}
\eeq
is equal to $\Z^d$. This is an  aperiodicity property. To see it is true, note that we can use the decay bound \eqref{eq:Kdeltadecay}, to conclude that for all sufficiently small $\de>0$, we have $T(0)>0$ and $T(\pm e_j)>0$ for $j=1,\ldots,d$.

\item The decay bound \eqref{eq:Kdeltadecay} also implies the summability of
\beq\label{eq:Uass2}
\begin{aligned}
&\sum_{x\in\Z^d}  |T(x)||x|^{2+m_d} <\infty,\qquad &d=3 \textnormal{ or } d\geq 5,\\
&\sum_{x\in\Z^d}  |T(x)||x|^{2+m_d} \ln |x| <\infty,\qquad &d=4,
\end{aligned}
\eeq
were $m_d=2d-3$, $ d=3,4 $  and $ m_{d}=d+1 $, $ d\ge 5 $, was defined in \eqref{eq:mddefn}. 

%Since the absolute convergence of a series implies its convergence, we may also conclude that
%\beq\label{eq:Uass2}
%\begin{aligned}
%&\sum_{x\in\Z^d} T(x)|x|^{2+m_d} <\infty,\qquad &d=3 \textnormal{ or } d\geq 5,\\
%&\sum_{x\in\Z^d} T(x)|x|^{2+m_d} \ln |x| <\infty,\qquad &d=4,
%\end{aligned}
%\eeq

\end{enumerate}

Together (i)-(iii) verify the assumptions of Theorem 2 in \cite{Uch98} with $m$ equal to $m_d$ and $T$ called $p$ there, with the exception of non-negativity. \medskip

\subsection{Verification of non-vanishing condition} We confirm that the fact that $T$ may be negative does not pose any problems in the proof. This step uses Proposition~\ref{prop:Qsymm} and the extension is applicable as long as $\mathbf Q$ is strictly positive semidefinite.

The proof of \cite[Theorem 2]{Uch98} is contained in Section 4 of that paper. 
The proof makes use of  general estimates on Fourier integrals taken from Sections 2 and 3 of \cite{Uch98} which do not depend on the non-negativity of $ p $. This concerns  Lemma~2.1, Lemma~3.1 and Corollary 3.1 from \cite{Uch98}. These are used in Section~4
together with the absolute summability \eqref{eq:Uass2} to control the error terms. The only step where the loss of non-negativity requires a short argument is the proof of the non-vanishing condition $c(\theta)^2+s(\theta)^2>0$ which is obtained on page 226 of \cite{Uch98} from positivity and aperiodicity. We now verify this condition to our context.
 
For $\theta\in\mathbb T^d$, we set
 $$
 c(\theta)=\sum_{x\in\Z^d} T(x)(1-\cos(\theta\cdot x)),\qquad  s(\theta)=\sum_{x\in\Z^d} T(x)\sin(\theta\cdot x).
 $$
 
 \be{lemma}\label{lm:cstheta} 
For sufficiently small $\de>0$, we have $$c(\theta)^2+s(\theta)^2>0 $$ for $\theta\in\mathbb T^d\setminus\{0\}$.
 \e{lemma}

 \be{proof}
 By Taylor expansion around $\theta=0$, we obtain
 $$
 \begin{aligned}
 c(\theta)=&\frac{1}{2}\sum_{x\in\Z^d} T(x) (\theta\cdot x)^2+\OO(\theta^4),\\
 s(\theta)=&\sum_{x\in\Z^d} T(x) (\theta\cdot x)+\OO(\theta^5),
 \end{aligned}
 $$
 where the error terms are controlled by the decay bound of $ \mathbf{K}^{\delta} $, \eqref{eq:Kdeltadecay} which  enters via the definition of $ T $ given in \eqref{eq:Tdefn}. By the definitions of $ \mathbf Q $,  \eqref{eq:Qdefn}, and $ T $, \eqref{eq:Tdefn}, we have $$  \mathbf{Q}=\mathrm{Hess}(m)(0)=-4d \mathrm{Hess}(\hat T)(0). $$ Using this for $ c $ and the zero mean property  \eqref{eq:Uass1} of $ T $, for $ s $ we obtain $$
 c(\theta)=\frac{1}{2d}\scp{\theta}{\mathbf{Q}\theta}+\OO(\theta^4),
 \qquad s(\theta)=\OO(\theta^5).
 $$
 Now Proposition~\ref{prop:Qsymm} says that for sufficiently small $\de>0$ the matrix $\mathbf Q$ is positive definite, i.e.,
  $$
 \scp{\theta}{\mathbf{Q}\theta}\geq (1-\de^2 C_d)\theta^2
 $$
 and therefore $c(\theta)^2+s(\theta)^2>0$ holds for all $|\theta|<r$ for some small $r>0$. 
 
 It remains to prove a lower bound over the set  $K=\setof{\theta\in\mathbb T^d}{|\theta|\geq r}$. To this end, let $c_0(\theta),s_0(\theta)$ denote the analogs of $c(\theta),s(\theta)$ with $\de=0$. Then we have $c_0(\theta)^2+s_0(\theta)^2>0$ for all $\theta\in K$ by the aperiodicity of the simple random walk. On the one hand, the continuous function $c_0(\theta)^2+s_0(\theta)^2$ takes its minimum on the compact set $K$; call it $\mu>0$. On the other hand, by \eqref{eq:Tdefn'} and \eqref{eq:Kdeltadecay}, we have
 $$
 \sup_{\theta\in K}|c(\theta)^2+s(\theta)^2-\l(c_0(\theta)^2+s_0(\theta)^2\r)|\leq \de^2 C_{d}.
 $$
 Thus, choosing $\de$ small enough that $\de^2 C_{d}\leq \mu/2$, we conclude that $c(\theta)^2+s(\theta)^2>0$ holds on $K$ as well. This proves Lemma \ref{lm:cstheta}.
 \e{proof}

\subsection{Conclusion} We are now ready to complete the proof of Theorem \ref{thm:GFasymptotic}. Thanks to (i)-(iii), Lemma \ref{lm:cstheta} and the paragraph preceding it, the proof of Theorem 2 from \cite{Uch98} extends to our situation and yields an asymptotic expansion similar to \eqref{eq:GFasymptotic}. Namely,  taking account of the rescaling by $4d$ that we introduced in \eqref{eq:mrewrite}, we have the asymptotic expansion
\beq\label{eq:Uchexpansion}
\mathcal G(x)=\frac{1}{4d}\frac{\kappa_d}{(\sigma')^2} |x'|^{2-d}
+\sum_{k=1}^{m_d} U_k\l(\frac{x'}{|x'|}\r) |x'|^{2-d-k} +o(|x'|^{2-d-m_d}), \quad \textnormal{as } |x|\to\infty,
\eeq
where $x'=\sigma'(\mathbf{Q}')^{-1/2} x$ and $\mathbf{Q}'$ is the $d\times d$ matrix generating the second-moment functional
$$
\scp{\theta}{\mathbf Q' \theta} = \sum_{x\in\Z^d} T(x) (x\cdot \theta)^2,\qquad \sigma'=(\det\mathbf Q')^{1/(2d)}.
$$

When we compare this with our claim  \eqref{eq:GFasymptotic}, we see that the latter features $\tilde x =\sigma\mathbf{Q}^{-1/2} x$ instead, with the matrix $\mathbf{Q}$ defined in \eqref{eq:Qdefn}. These are related via
\beq\label{eq:Q'Qconnect}
\mathbf{Q}'=\frac{1}{4d}\mathbf Q
\eeq
To see this, we use the Fourier representation and recall \eqref{eq:mrewrite} and \eqref{eq:mthetadefn} to find
$$
 \begin{aligned}
\mathbf{Q}_{i,j}'
=\sum_{x\in\Z^d} T(x) x_i x_j
=- \l(\frac{\partial^2}{\partial \theta_i\partial \theta_j}\hat T\r)(0)
=\frac{1}{4d}\l(\frac{\partial^2}{\partial \theta_i\partial \theta_j}m\r)(0)
= \frac{1}{4d} \mathbf{Q}_{i,j}
 \end{aligned}
$$
for all $i,j\in\{1,\ldots,d\}$.

Finally, we employ the identity \eqref{eq:Q'Qconnect} and its consequence $\sigma'=(4d)^{-1/2}\sigma$ in \eqref{eq:Uchexpansion}. For the leading term, we note that $|x'|=|x|$ and $\frac{1}{4d(\sigma')^2}=\frac{1}{\sigma^2}$. For the subleading term, we note that $\frac{x'}{|x'|}=\frac{\tilde x}{|\tilde x|}$. Absorbing the factors of $(4d)^{d+k-2}$ into the $U_1,\ldots,U_m$ then yields \eqref{eq:GFasymptotic}. This proves Theorem \ref{thm:GFasymptotic}.
 \qed\\

\section{Proof of Corollary \ref{cor:nablaGFasymptotic} and Formula \eqref{eq:nablaGFasymptotic}}
\label{sect:maincorpf}

\be{proof}[Proof of Corollary \ref{cor:nablaGFasymptotic}] 
Let $1\leq j\leq d$, $s\in\{ \pm 1\}$ and let $\alpha\in \N_0^d$ be a multi-index with $|\alpha|\leq m_d$. We apply Theorem \ref{thm:GFasymptotic}. Regarding the $U_k$ terms, \eqref{eq:Uorder} shows that 
$$
\nabla^\alpha \l(\sum_{k=m_d-|\al|+1}^{m_d} U_k\l(\frac{\tilde x}{|\tilde x|}\r) |\tilde x|^{2-d-k}\r)\in O(|\tilde x|^{1-d-m_d})\subseteq o(|\tilde x|^{2-d-m_d}).
$$
Regarding the $o(|\tilde x|^{2-d-m_d})$ error term appearing in Theorem \ref{thm:GFasymptotic}, we note that 
for $f\in o(|\tilde x|^{2-d-m_d})$, the triangle inequality implies $\nabla^\alpha f([x,x+se_j])
\in o(|\tilde x|^{2-d-m_d})$, so the error does not get worse under discrete differentiation. This proves Corollary~\ref{cor:nablaGFasymptotic}.
\e{proof}

\be{proof}[Proof of Formula \eqref{eq:nablaGFasymptotic}]
We recall the notation \eqref{eq:tildexdefn}, i.e., $\tilde x=\sigma \mathbf{Q}^{-1/2} x$. We first use Corollary \ref{cor:nablaGFasymptotic} and the bound \eqref{eq:Uorder} for all $k\geq 1$ to find
\begin{align*}
\nabla_j \mathcal G(x)=\frac{\kappa_d}{\sigma^2} ( |\tilde x+\tilde e_j|^{2-d}- |\tilde x|^{2-d}) +\OO(|\tilde x|^{-d})
\end{align*}
To compute the leading term, we expand $(1+y)^q=1+qy+\OO(y^2)$ to obtain
$$
\begin{aligned}
 |\tilde x+\tilde e_j|^{2-d}- |\tilde x|^{2-d}
 &=|\tilde x|^{2-d}\l(\l(1+\scpp{\frac{\tilde x}{|\tilde x|^2}}{\tilde e_j}+\frac{|\tilde e_j|^2}{|\tilde x|^2}\r)^{\frac{2-d}{2}}-1\r)\\
 &=|\tilde x|^{1-d} \frac{2-d}{2} \scpp{\frac{\tilde x}{|\tilde x|}}{\tilde e_j}+\OO(|\tilde x|^{-d}),\qquad \textnormal{as } |x|\to\infty,
\end{aligned}
$$
where we also made use of the equivalence of  the norms $|x|  $ and $ |\tilde x| $, cf.\ \eqref{eqn:xtilde}.
\e{proof}

\section*{Acknowledgments} The authors are grateful to Scott Armstrong and Mitia Duerinckx for useful remarks. MK acknowledges the financial support of the German Science Foundation. 

\begin{appendix}
\section{Direct argument for the leading order in Theorem \ref{thm:GFasymptotic}}\label{sect:direct}
In this appendix, we give a self-contained proof of the lowest order asymptotic in Theorem \ref{thm:GFasymptotic}, i.e.,
\beq\label{eq:direct}
\mathcal G(x)=\frac{\kappa_d}{\sigma^2}|\tilde x|^{2-d} +\OO(|x|^{1-d}),\qquad \textnormal{as } |x|\to\infty,
\eeq
where again  $\tilde x=\sigma \mathbf{Q}^{-1/2} x$.
The approach is to use a Taylor expansion around the origin in Fourier space which is justified by Theorem \ref{thm:bkl} and controlled by adapting the dyadic pigeonholing from \cite[Appendix A]{KL}.

\textit{Proof of \eqref{eq:direct}.}
We recall Definition \eqref{eq:mthetadefn} of $m(\theta)$ and the fact that $\widehat{K^\de_{j,k}}\in C^{2d-1}(\mathbb T^d)$ for all $j,k\in \{1,\ldots,d\}$. We first isolate the lowest, quadratic order of $m(\theta)$ by setting
$$
m(\theta)=m_0(\theta)+\tilde m(\theta),
\qquad \textnormal{where }
m_0(\theta)=\scpp{\theta}{\mathbf{Q}\theta}.
$$
Observe that $\tilde m\in C^{2d-1}(\mathbb T^d)$ satisfies
\beq\label{eq:tildembounds}
|D^\al \tilde m(\theta)|\leq C_{d}|\theta|^{3-|\al|},\qquad  \al\in \mathbb N_0^d \textnormal{ with } |\al|\leq 2d-1.
\eeq
We can decompose
\beq\label{eq:1/mdecomp}
\frac{1}{m(\theta)} =\frac{1}{m_0(\theta)} -\frac{\tilde m(\theta)}{m_0(\theta)m(\theta)}.
\eeq
The next lemma then implies \eqref{eq:direct}.

%Combining \eqref{eq:1/mdecomp} with \eqref{eq:GavgFourierrep}, \eqref{eq:G_lt} and \eqref{eq:G_hot} yields \eqref{eq:direct}.

\be{lemma}\label{lm:GFsubleading}
For all $\de\geq 0$ sufficiently small, we have
\beq\label{eq:G_hot}
\int_{\mathbb T^d} e^{ix\cdot \theta} \frac{\tilde m(\theta)}{m_0(\theta)m(\theta)}\frac{\mathrm{d}^d \theta}{(2\pi)^d}=\OO(|x|^{1-d}),\qquad \textnormal{as } |x|\to\infty.
\eeq
\beq\label{eq:G_lt}
\int_{\mathbb T^d} e^{ix\cdot \theta} \frac{1}{m_0(\theta)} \frac{\mathrm{d}^d \theta}{(2\pi)^d}
=\frac{\kappa_d}{\sigma^2} |\tilde x|^{2-d} +\OO(|x|^{1-d}),\qquad \textnormal{as } |x|\to\infty.
\eeq
\e{lemma}

\begin{proof}
We loosely follow Appendix A in \cite{KL} where a similar problem is treated and start with the proof of \eqref{eq:G_hot}. We define the function $ F $ on $ \mathbb{T}^{d} $ by
$$
F(\theta)=\frac{\tilde m(\theta)}{m_0(\theta)m(\theta)}.
$$
and note that, due to \eqref{eq:tildembounds} and the quadratic vanishing order of $m(\theta)$ and $m_0(\theta)$ at the origin, we have
\beq\label{eq:DalFestimate}
|F(\theta)|\leq C_d |\theta|^{-1}.
\eeq
%In essence, we aim to use integration by parts to transfer $d-1$ derivatives onto $F$. This leads to a singular Fourier integral, which we analyze through the following dyadic scale decomposition of the frequencies. 
Let $\varphi:[0,\infty)\to\R$ be a smooth cutoff function with $\varphi=1$ on $[0,2\pi]$ which is supported on $[0,4\pi]$. Define $\psi(r):=\varphi(r)-\varphi(2r)$ and $\psi_l(r)=\psi(2^l r)$ for all $l\geq 1$ and $ r\ge 0 $. Note that this defines a partition of unity $\sum_{l\geq 0}\psi_l(r)=1$ for all $r\neq 0$. We decompose
\beq\label{eq:fldecompose}
\int_{\mathbb T^d} e^{ix\cdot \theta} F(\theta)\frac{\mathrm{d}^d \theta}{(2\pi)^d}
=\sum_{l\geq 0} f_l(x).
\eeq
where we rescaled and introduced
$$
f_l(x)=2^{-ld}\int_{\mathbb T^d} e^{i2^{-l}x\cdot \theta} F_l(\theta)\frac{\mathrm{d}^d \theta}{(2\pi)^d},\qquad F_l(\theta)=\psi(|\theta|) F(2^{-l}\theta).
$$
Note that \eqref{eq:DalFestimate} implies $|F_l(\theta)|\leq C_d 2^{l} |\theta|^{-1}$, so we can use the triangle inequality to obtain
\beq\label{eq:flestimate1}
|f_l(x)|\leq C_d 2^{-l(d-1)}.
\eeq
This is useful for for $|x|2^{-l}\leq 1$, while for $|x|2^{-l}\geq 1$ it can be improved to
\beq\label{eq:flestimate2}
|f_l(x)|\leq C_d \frac{2^{-l(d-1)}}{(|x|2^{-l})^{d}}.
\eeq
To prove \eqref{eq:flestimate2}, we assume without loss of generality that $|x_1|=\max_{1\leq j\leq d}|x_j|$ and use $d$-fold integration by parts to write
\beq\label{eq:flibp}
f_l(x) =i^d \frac{2^{-ld}}{(x_1 2^{-l})^d}
\int_{\mathbb T^d} e^{i2^{-l}x\cdot \theta} \partial_{\theta_1}^d F_l(\theta)\frac{\mathrm{d}^d \theta}{(2\pi)^d}.
\eeq
Next, we observe that \eqref{eq:tildembounds} and the quadratic behavior of both $m_0$ and $m$ at the origin yield that $\|\partial_{\theta_1}^d F_l\|_\infty\leq C_d 2^{-l(d-1)}$; compare Lemma A.1 in \cite{KL}. Applying this bound to \eqref{eq:flibp} yields \eqref{eq:flestimate2}.

We use \eqref{eq:flestimate1} and \eqref{eq:flestimate2} and bound the resulting geometric series to find
$$
\begin{aligned}
\sum_{l\geq 0} |f_l(x)|
\leq C_d \sum_{l\geq \log_2 |x|} 2^{-l(d-1)}+C_d \sum_{0\leq l\leq \log_2 |x|}\frac{2^{-l(d-1)}}{(|x|2^{-l})^{d}}
\leq C_d |x|^{1-d}.
\end{aligned}
$$
In view of \eqref{eq:fldecompose}, this proves \eqref{eq:G_hot}.

Next, we turn to the proof of \eqref{eq:G_lt}. By Proposition~\ref{prop:Qsymm}, the matrix $\mathbf{Q}$ is symmetric and positive definite. Hence, by a change of variables,
$$
\begin{aligned}
\int_{\mathbb T^d} e^{ix\cdot \theta} \frac{1}{m_0(\theta)} \frac{\mathrm{d}^d \theta}{(2\pi)^d}
%=&\int_{\mathbb T^d} e^{ix\cdot \theta} \frac{1}{\scp{\mathbf{Q}^{1/2}\theta}{\mathbf{Q}^{1/2}\theta}} \frac{\mathrm{d}^d \theta}{(2\pi)^d}\\
=&\sigma^{-d}
\int_{\mathbb T^d} e^{i (\mathbf{Q}^{-1/2}x)\cdot \theta} \frac{1}{|\theta|^2} \frac{\mathrm{d}^d \theta}{(2\pi)^d}.
\end{aligned}
$$
For the volume element, we used that $\det(\mathbf{Q}^{1/2})= (\det\mathbf{Q})^{1/2}=\sigma^d$ since $\mathbf{Q}$ is symmetric.
%Next we leverage the homogeneity $|\theta|^{-2}$ of the singularity by extending the integral from the torus $\mathbb T^d$ to the whole Euclidean space $\R^d$.
By applying \eqref{eq:G_hot} with $\de=0$, we obtain
$$
\begin{aligned}
\int_{\mathbb T^d} e^{i \mathbf{Q}^{-1/2}x\cdot \theta} \frac{1}{|\theta|^2} \frac{\mathrm{d}^d \theta}{(2\pi)^d}
=\int_{\mathbb T^d} e^{i \mathbf{Q}^{-1/2}x\cdot \theta} \frac{1}{2\sum_{j=1}^d (1-\cos\theta_j)} \frac{\mathrm{d}^d \theta}{(2\pi)^d}
+\OO(|x|^{1-d}).
\end{aligned}
$$
We recognize the first integral as the Green's function of the free Laplacian on $\Z^d$ evaluated at the point $\mathbf Q^{-1/2}x$. The standard asymptotic formula for the free Laplacian gives
$$
\int_{\mathbb T^d} e^{ix\cdot \theta} \frac{1}{m_0(\theta)} \frac{\mathrm{d}^d \theta}{(2\pi)^d}
=\frac{\kappa_d}{\sigma^d} |\mathbf{Q}^{-1/2}x|^{2-d} +\OO(|x|^{1-d})
=\frac{\kappa_d}{\sigma^2} |\tilde x|^{2-d} +\OO(|x|^{1-d}).
$$
This proves \eqref{eq:G_lt} and thus Lemma \ref{lm:GFsubleading}.
\end{proof}

\end{appendix}

\bibliographystyle{amsplain}

\end{document}